\documentclass[journal,twocolumn,letterpaper]{IEEEJERM}

\ifCLASSINFOpdf
\else
\fi
\usepackage{cite}
\usepackage{amsmath,amssymb,amsfonts}
\usepackage{algorithmic}
\usepackage{graphicx}
\usepackage{textcomp}
\usepackage{subfigure}

\begin{document}
\title{A High-Order Discretization Scheme for Surface Integral Equations for Analyzing the Electroencephalography Forward Problem}

\author{Rui~Chen,~\IEEEmembership{Member,~IEEE,} Viviana~Giunzioni,~\IEEEmembership{Member,~IEEE,}
Adrien~Merlini,~\IEEEmembership{Senior Member,~IEEE,}
and~Francesco~P. ~Andriulli,~\IEEEmembership{Fellow,~IEEE}
\thanks{This work was supported in part by the National Natural Science Foundation of China (NSFC) under Grant 62201264 and Grant 62331016; in part by the Fundamental Research Funds for the Central Universities under Grant 30924010207; in part by the Fund Program for the Scientific Activities of Selected Returned Overseas Professionals in Shanxi Province under Grant 20240063; in part by European Innovation Council (EIC) under Grant 101046748; in part by European Union-Next Generation EU under Grant CUP E13C22000990001; and in part by the Labex CominLabs Excellence Laboratory under Grant ANR-10-LABX-07-01.  \emph{(Corresponding author: Rui Chen.)}}
\thanks{Rui Chen is with the School of Microelectronics, Nanjing University of Science and
Technology, Nanjing 210094, China (e-mail: rui.chen@njust.edu.cn). }
\thanks{Viviana Giunzioni and Francesco P. Andriulli are with the Department of Electronics and Telecommunications, Politecnico di Torino, 10129 Turin, Italy (e-mail: viviana.giunzioni@polito.it; francesco.andriulli@polito.it).}
\thanks{Adrien Merlini is with the Department of Microwave, IMT Atlantique, 29238 Brest, France (e-mail: adrien.merlini@imt-atlantique.fr).}
}



\maketitle

\begin{abstract}

A Nyström-based high-order (HO) discretization scheme for surface integral equations (SIEs) for analyzing the electroencephalography (EEG) forward problem is proposed in this work. We use HO surface elements and interpolation functions for the discretization of the interfaces of the head volume and the unknowns on the elements, respectively. The advantage of this work over existing isoparametric HO discretization schemes resides in the fact that the interpolation points are different from the mesh nodes, allowing for  the flexible manipulation of the order of the basis functions without regenerating the mesh of the interfaces. Moreover, the interpolation points are chosen from the quadrature rules with the same number of points on the elements simplifying the numerical computation of the surface integrals  for the far-interaction case. In this contribution, we extend the implementation of the HO discretization scheme to the double-layer and the adjoint double-layer formulations, as well as to the isolated-skull-approach for the double-layer formulation and to the indirect adjoint double-layer formulation, employed to improve the solution accuracy in case of high conductivity contrast models, which requires the development of different techniques for the singularity treatment. Numerical experiments are presented to demonstrate the accuracy, flexibility, and efficiency of the proposed scheme for the four SIEs for analyzing the EEG forward problem.

\end{abstract}




%

\section{Introduction}
\label{sec:introduction}
\IEEEPARstart{E}{lectroencephalography} (EEG) is one of the widely used non-invasive neurophysiological techniques for functional neuroimaging that aims to determine location, intensity, and orientation of the activity of the neural source in the head volume \cite{Cabeza}. Although EEG-based neuroimaging has good temporal resolution, its spatial resolution needs further improvement \cite{Acar2016}.

EEG-based imaging techniques characterize the neural activity in the head volume from the  electric potential directly measured at the electrodes placed on the surface of the scalp, which is termed as ``EEG inverse problem'' \cite{Phillips}. To analyze the EEG inverse problem, the solution of the EEG forward problem, i.e., the evaluation of the electric potential on the scalp surface with the given hypothetical neural source and electrical property of the biological tissues in the head volume, is required \cite{Marqui}. Since the accuracy of the analysis of the EEG inverse problem depends on that of the solution of the EEG forward problem \cite{Acar2013}, we could improve the spatial resolution of EEG via the later.

In the last decades, many numerical methods have been developed for the solution of the EEG forward problem. Two kinds of the commonly used methods are differential equation- and surface integral equation (SIE)-based  methods \cite{Awada, Piastra, Vorwerk, Yavich, He, Hallez, Kybic, Munck, Schlitt, Meijs, Hamalainen, Giunzioni,Rahmouni, Monin, Frijns, Gençer1999, Gençer2005}. Compared with the differential equation-based methods \cite{Awada, Piastra, Vorwerk, Yavich}, the SIE-based methods \cite{He, Hallez, Kybic, Munck, Schlitt, Meijs, Hamalainen, Frijns, Gençer1999, Gençer2005,Giunzioni,Rahmouni, Monin} remove the need of (artificial) absorbing boundary conditions and use surface elements instead of volumetric elements for the object discretization, resulting in fewer unknowns \cite{He}.

The SIE-based methods solve SIEs enforced on the interfaces between the neighboring biological tissue layers in the head volume for the unknowns on the interfaces \cite{Hallez}. To improve the accuracy of the solution of the EEG forward problem, one could improve that of the discretization of the interfaces and unknowns \cite{Kybic}. On one hand, traditionally, low-order (LO) surface elements (e.g., planar triangular patch) are usually used for the  discretization of the  interfaces  \cite{ He, Hallez, Kybic, Munck, Schlitt, Meijs, Hamalainen, Giunzioni,Rahmouni, Monin}. However, as the interfaces are not regular, the use of LO surface elements introduces non-negligible geometrical  error. Instead, the use high-order (HO) surface elements (e.g., curved triangular patch) improves the modeling accuracy \cite{Frijns, Gençer1999, Gençer2005}. On the other hand, traditionally, LO basis functions are  usually used for the discretization of the unknowns \cite{ He, Hallez, Kybic, Munck, Schlitt, Meijs, Hamalainen, Giunzioni,Rahmouni, Monin}. However, the approximation accuracy is limited using LO basis functions. Instead, one could use HO basis functions to improve the approximation accuracy of the unknowns \cite{Frijns, Gençer1999, Gençer2005}. 

In the past years, several works have used the HO surface elements and HO basis functions simultaneously for the discretization of SIEs for analyzing the EEG forward problem \cite{Frijns, Gençer1999, Gençer2005}. \cite{Frijns} used quadratically curved triangular elements for the discretization of the interfaces and quadratic interpolation functions associated with the interpolation points for the discretization of the electric potential.  \cite{Gençer1999} and \cite{Gençer2005} used the HO basis functions for the discretization of the unknown on the planar, quadratic, and cubic triangular elements.  \cite{Frijns,Gençer1999, Gençer2005} used the isoparametric HO discretization scheme for SIEs, i.e., the same interpolation functions are used for the definition of the surface elements and of the unknown on the elements and the mesh nodes are overlapped with the interpolation points for each element. However, it limits the flexibility of the HO discretization of SIEs. For instance, to use the higher-order basis functions, one has to first regenerate the mesh of the interfaces with the higher-order surface elements. 

To address this issue, we propose a Nyström-based HO discretization scheme for SIEs for analyzing the EEG forward problem in this work. The unknowns on the HO surface elements are approximated using the HO interpolation functions associated with the interpolation points. The novelty of this work over the existing works using the isoparametric HO discretization scheme is ``threefold'': 
\begin{enumerate}
\item The interpolation points are located on the elements instead of the boundary of the elements. Since the interpolation points are independent of the mesh nodes, we could flexibly manipulate the order of the basis functions without regenerating the mesh of the interfaces.
\item The interpolation points are chosen from the quadrature rules with the same number of points on the elements. The numerical computation of the surface integrals can be simplified for the far-interaction case. 
\item We extend the implementation of the HO discretization scheme to the double-layer approach and the adjoint double-layer approach, as well as to the isolated skull approach for the double-layer formulation and to the indirect adjoint double-layer approach for  analyzing the EEG forward problem, which needs the development of different  techniques for the singularity treatment. 
\end{enumerate}
Note that, even though the Nyström-based methods have been reported in electromagnetics \cite{Kang, Chen1, Chen2} and acoustics \cite{Canino, Chen3, Chen4}, to the best of the authors' knowledge, the Nyström-based HO discretization scheme has never been used for SIEs for analyzing the EEG forward problem. A preliminary conference version of this work is presented in \cite{Chen5} while this paper  presents more implementation details.


\section{Formulation}
\label{sec:formulation}

\subsection{SIEs for the EEG Forward Problem}
\label{subsec:SIEs}
\begin{figure}[!t]
\centering
\includegraphics[width=0.9\columnwidth]{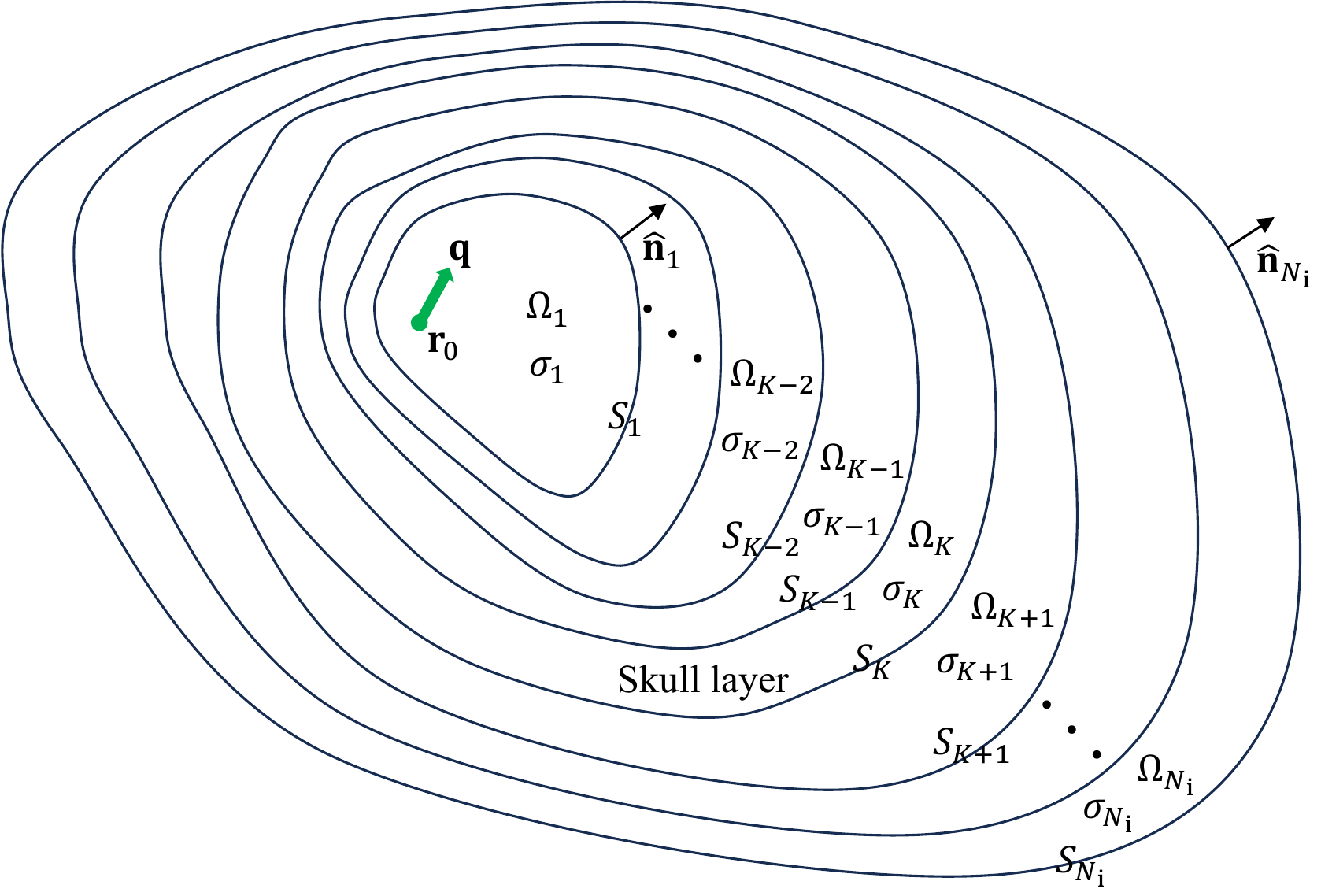}
\caption{The head volume model. }
\label{fig_1}
\end{figure}

When using SIEs for analyzing the EEG forward problem, the head volume $\Omega$ is usually modeled as a piecewise, homogeneous, and multilayered conductor (as seen in Fig. \ref{fig_1}). The layers (representing the biological tissues) from the inner to the outer are denoted as $\Omega_1$, $\Omega_2$, ..., $\Omega_{N_{\mathrm{i}}}$, with the constant conductivities $\sigma_1$, $\sigma_2$, ..., $\sigma_{N_{\mathrm{i}}}$, respectively. The skull layer is located in $\Omega_K$ between $S_{K-1}$ and $S_K$ with the conductivity $\sigma_K$. The interfaces between the neighboring layers are $S_1$, $S_2$, ..., $S_{N_{\mathrm{i}}}$ from the inner to the outer, where the outward pointing unit normal
vectors $\hat{\mathbf{n}}_1$, $\hat{\mathbf{n}}_2$, ..., $\hat{\mathbf{n}}_{N_{\mathrm{i}}}$ are defined. Note that we assume that the conductivity of the homogeneous background medium outside $\Omega$ is zero in this work. 

In addition, the neural source in the head volume is often modeled as an electric current dipole for SIEs for the EEG forward problem \cite{Kybic}. The electric potential  generated by the dipolar source in the  background region is 
\begin{align}
V^{\mathrm{inc}}(\mathbf{r})=\frac{\mathbf{q} \cdot (\mathbf{r}-\mathbf{r}_0)}{4 \pi |\mathbf{r}-\mathbf{r}_0|^3}
\label{eq_1}
\end{align}
where $\mathbf{q}$ is the dipolar moment, $\mathbf{r}$ is the observation point, and the dipole position is assumed to be $\mathbf{r}_0 \in \Omega_1$ in this work.

In the following, four approaches are introduced for computing the electric potential $V(\mathbf{r})$ for the EEG forward problem.

\subsubsection{Double-Layer Approach}
SIE in the double-layer (DL) formulation is expressed as \cite{Kybic} 
\begin{align}
\nonumber V^{\mathrm{inc}}&(\mathbf{r})=\frac{\sigma_j+\sigma_{j+1}}{2}V(\mathbf{r})
-\sum_{i=1}^{N_{\mathrm{i}}}(\sigma_{i+1}-\sigma_i) \times \\      & \mathrm{P. V.} \int_{S_i}\partial_{\mathbf{n}'}G(\mathbf{r},\mathbf{r}')V(\mathbf{r}')d\mathbf{r}' , \ \ \  \mathbf{r} \in S_j,\,\,j = 1,\dots,N_{\mathrm{i}}
\label{eq_2}
\end{align}
where $\mathbf{r}'$ is the source point, $\partial_{\mathbf{n}'}=  \hat{\mathbf{n}}'(\mathbf{r}') \cdot \nabla'$, $\hat{\mathbf{n}}'(\mathbf{r}')$ is the outward pointing unit normal vector at $\mathbf{r}'$, $G(\mathbf{r},\mathbf{r}')=1/(4\pi |\mathbf{r}-\mathbf{r}' |)$ is the Green function, and $V (\mathbf{r})$ and $V (\mathbf{r}')$ represent the electric potential at $\mathbf{r}$ and $\mathbf{r}'$, respectively. Here, ``P.V.'' denotes the ``principal value''.

\subsubsection{Adjoint Double-Layer Approach}
SIE in the adjoint double-layer (ADL) formulation is expressed as \cite{Kybic}
\begin{align}
\nonumber &\partial_{\mathbf{n}}\left[\frac{V^{\mathrm{inc}}(\mathbf{r})}{\sigma_1}\right]=\frac{\sigma_j+\sigma_{j+1}}{2(\sigma_{j+1}-\sigma_j)}\xi(\mathbf{r})\\
&-\sum_{i=1}^{N_{\mathrm{i}}} \mathrm{P. V.} \int_{S_i}\partial_{\mathbf{n}}G(\mathbf{r},\mathbf{r}')\xi(\mathbf{r}')d\mathbf{r}', \ \ \mathbf{r} \in S_j,\,\,j = 1,\dots,N_{\mathrm{i}}
\label{eq_3}
\end{align}
where  $\partial_{\mathbf{n}}=\hat{\mathbf{n}}(\mathbf{r}) \cdot  \nabla $, and $\xi (\mathbf{r})$ and $\xi (\mathbf{r}')$ represent the jump of $\partial_{\mathbf{n}} V (\mathbf{r})$ and $\partial_{\mathbf{n}} V (\mathbf{r}')$ along the direction of the outward pointing unit normal vectors at $\mathbf{r}$ and $\mathbf{r}'$, respectively.

After $\xi (\mathbf{r})$ is obtained by solving \eqref{eq_3}, $V (\mathbf{r})$ can be calculated as a sum of an electric potential term for excitation and single-layer surface integral operators for  correction as follows   \cite{Kybic}
\begin{align}
V(\mathbf{r})&=\frac{V^{\mathrm{inc}}(\mathbf{r})}{\sigma_1}\nonumber\\&+\sum_{i=1}^{N_{\mathrm{i}}}\int_{S_i}G(\mathbf{r},\mathbf{r}')\xi(\mathbf{r}')d\mathbf{r}', \ \ \mathbf{r} \in S_j,\,\,j = 1,\dots,N_{\mathrm{i}}.
\label{eq_4}
\end{align}

\subsubsection{Isolated Skull Approach for the Double-Layer Approach}
The isolated skull approach (ISA) is originally developed to address the numerical inaccuracies of the DL formulation due to the low conductivity of one layer of the head volume conductor model \cite{Hamalainen}. ISA decomposes $V (\mathbf{r})$  as 
\begin{align}
V (\mathbf{r})=
\begin{cases}
V_{\mathrm{ISA}} (\mathbf{r})+V_{\mathrm{corr}} (\mathbf{r}),  \ \ \mathrm{if} \  \mathbf{r} \in S_1,...,S_{K-1}\\
V_{\mathrm{corr}} (\mathbf{r}), \ \  \mathrm{if} \ \mathbf{r} \in S_K,S_{K+1},...,S_{N_{\mathrm{i}}}
\end{cases}
\label{eq_5}
\end{align}
where $V_{\mathrm{ISA}} (\mathbf{r})$ is the solution of an isolated model considering only the tissues under the skull layer with the surfaces $S_1, S_2, ..., S_{K-1}$ and $V_{\mathrm{corr}} (\mathbf{r})$ is the correction term. Note that, $V_{\mathrm{ISA}} (\mathbf{r})=0$ at $\mathbf{r} \in  S_K, S_{K+1}, ..., S_{N_{\mathrm{i}}}$ for ISA.

SIE in the unknown of $V_{\mathrm{ISA}} (\mathbf{r})$ is expressed as
\begin{align}
\nonumber &V^{\mathrm{inc}}(\mathbf{r})=\frac{\sigma_j+\sigma_{j+1}}{2}V_{\mathrm{ISA}}(\mathbf{r})-\sum_{i=1}^{K-2}(\sigma_{i+1}-\sigma_i) \times \\ \nonumber &\mathrm{P. V.} \int_{S_i}\partial_{\mathbf{n}'}G(\mathbf{r},\mathbf{r}')V_{\mathrm{ISA}}(\mathbf{r}')d\mathbf{r}'+\sigma_{K-1} \times \\
&\mathrm{P. V.} \int_{S_{K-1}}\partial_{\mathbf{n}'}G(\mathbf{r},\mathbf{r}')V_{\mathrm{ISA}}(\mathbf{r}')d\mathbf{r}', \ \mathbf{r} \in S_1, ..., S_{K-1}
\label{eq_6}
\end{align}
and SIE in the unknown of $V_{\mathrm{corr}} (\mathbf{r})$ is expressed as
\begin{align}
\nonumber &V_{\mathrm{corr}}^{\mathrm{inc}}(\mathbf{r})=\frac{\sigma_j+\sigma_{j+1}}{2}V_{\mathrm{corr}}(\mathbf{r})
-\sum_{i=1}^{N_{\mathrm{i}}}(\sigma_{i+1}-\sigma_i) \times \\  & \mathrm{P. V.} \int_{S_i}\partial_{\mathbf{n}'}G(\mathbf{r},\mathbf{r}')V_{\mathrm{corr}}(\mathbf{r}')d\mathbf{r}' , \  \mathbf{r} \in S_1, S_2, ..., S_{N_{\mathrm{i}}}
\label{eq_7}
\end{align}
where the expression of $V_{\mathrm{corr}}^{\mathrm{inc}}(\mathbf{r})$ is given by
\begin{align}
\nonumber &V_{\mathrm{corr}}^{\mathrm{inc}}(\mathbf{r})=\\
&\begin{cases}
    \frac{\sigma_K}{\sigma_{K-1}} V^{\mathrm{inc}}(\mathbf{r})+
    \sum_{i=1}^{K-2} \frac{\sigma_K (\sigma_{i+1}-\sigma_i)}{\sigma_{K-1}} \times \\ \ \ \ \   \int_{S_i}\partial_{\mathbf{n}'}G(\mathbf{r},\mathbf{r}')V_{\mathrm{ISA}}(\mathbf{r}')d\mathbf{r}' , \ \ \mathrm{if} \ \mathbf{r} \in S_K, S_{K+1}, ..., S_{N_{\mathrm{i}}} \\
    \frac{\sigma_K}{\sigma_{K-1}} V^{\mathrm{inc}}(\mathbf{r})-\sigma_K V_{\mathrm{ISA}}(\mathbf{r})+\sum_{i=1}^{K-2} \frac{\sigma_K (\sigma_{i+1}-\sigma_i)}{\sigma_{K-1}} \times \\ \ \ \ \  \int_{S_i}\partial_{\mathbf{n}'}G(\mathbf{r},\mathbf{r}')V_{\mathrm{ISA}}(\mathbf{r}')d\mathbf{r}' , \ \ \mathrm{if} \ \mathbf{r} \in  S_{K-1} \\
    \frac{\sigma_K}{\sigma_{K-1}} V^{\mathrm{inc}}(\mathbf{r})-\frac{\sigma_K (\sigma_j+\sigma_{j+1})}{2 \sigma_{K-1}} V_{\mathrm{ISA}}(\mathbf{r})+\sum_{i=1}^{K-2} \frac{\sigma_K (\sigma_{i+1}-\sigma_i)}{\sigma_{K-1}}  \\ 
    \ \ \ \   \times \mathrm{P. V.} \int_{S_i}\partial_{\mathbf{n}'}G(\mathbf{r},\mathbf{r}')V_{\mathrm{ISA}}(\mathbf{r}')d\mathbf{r}' , \ \ \mathrm{if} \ \mathbf{r} \in  S_1, ..., S_{K-2} 
\end{cases}
\label{eq_8}
\end{align}

\subsubsection{Indirect Adjoint Double-Layer Approach}
The indirect adjoint double-layer (IADL) formulation is originally developed to address the high conductivity ratio issue of the ADL formulation \cite{Rahmouni}. IADL writes $V(\mathbf{r})$, $\mathbf{r} \in \Omega_{N_{\mathrm{i}}}$ using auxiliary sources $J_{N_{\mathrm{i}}-1}(\mathbf{r})$, $\mathbf{r} \in S_{N_{\mathrm{i}}-1}$ and $J_{N_{\mathrm{i}}}(\mathbf{r})$, $\mathbf{r} \in S_{N_{\mathrm{i}}}$ as
\begin{align}
\nonumber    V(\mathbf{r})&=\int_{S_{N_{\mathrm{i}}-1}} G(\mathbf{r},\mathbf{r}') J_{N_{\mathrm{i}}-1}(\mathbf{r}') d
 \mathbf{r}' \\
 &+ \int_{S_{N_{\mathrm{i}}}} G(\mathbf{r},\mathbf{r}') J_{N_{\mathrm{i}}}(\mathbf{r}') d
 \mathbf{r}', \ \ \  \mathbf{r} \in \Omega_{N_{\mathrm{i}}}.
\label{eq_9}
\end{align}

Additionally, SIEs in unknowns of $J_{N_{\mathrm{i}}-1}(\mathbf{r})$ and $J_{N_{\mathrm{i}}}(\mathbf{r})$ are 
\begin{align}
\nonumber  &  \frac{\sigma_{N_{\mathrm{i}}-1}}{\sigma_{N_{\mathrm{i}}}-\sigma_{N_{\mathrm{i}}-1}} \xi (\mathbf{r}) = \\
\nonumber &- \frac{1}{2} J_{N_{\mathrm{i}}-1}(\mathbf{r})
+ \mathrm{P. V.} \int_{S_{N_{\mathrm{i}}-1}} \partial_{\mathbf{n}}G(\mathbf{r},\mathbf{r}') J_{N_{\mathrm{i}}-1}(\mathbf{r}') d \mathbf{r}'\\
    &+\int_{S_{N_{\mathrm{i}}}} \partial_{\mathbf{n}}G(\mathbf{r},\mathbf{r}') J_{N_{\mathrm{i}}}(\mathbf{r}') d \mathbf{r}', \ \ \  \mathbf{r} \in S_{N_{\mathrm{i}}-1}^{+} \label{eq_10}\\
\nonumber    &  0 =  \frac{1}{2} J_{N_{\mathrm{i}}}(\mathbf{r})+\int_{S_{N_{\mathrm{i}}-1}} \partial_{\mathbf{n}}G(\mathbf{r},\mathbf{r}') J_{N_{\mathrm{i}}-1}(\mathbf{r}') d \mathbf{r}' \\
      &+\mathrm{P. V.} \int_{S_{N_{\mathrm{i}}}} \partial_{\mathbf{n}}G(\mathbf{r},\mathbf{r}') J_{N_{\mathrm{i}}}(\mathbf{r}') d \mathbf{r}', \ \ \  \mathbf{r} \in S_{N_{\mathrm{i}}}^{-}   \label{eq_11}
\end{align}
where $\xi (\mathbf{r})$ is the solution of \eqref{eq_3}. After the solution of \eqref{eq_10}-\eqref{eq_11}, $V(\mathbf{r})$, $\mathbf{r} \in \Omega_{N_{\mathrm{i}}}$ could be obtained using \eqref{eq_9}.

\subsection{Discretization of SIEs}
\label{subsec:discretization}

To numerically solve SIEs  in Section  \ref{subsec:SIEs}, we use  the Nyström method \cite{Kang} in this work. The basic  implementation steps of the Nyström method are introduced as follows.

\begin{figure}[!t]
\centering
\includegraphics[width=0.9\columnwidth]{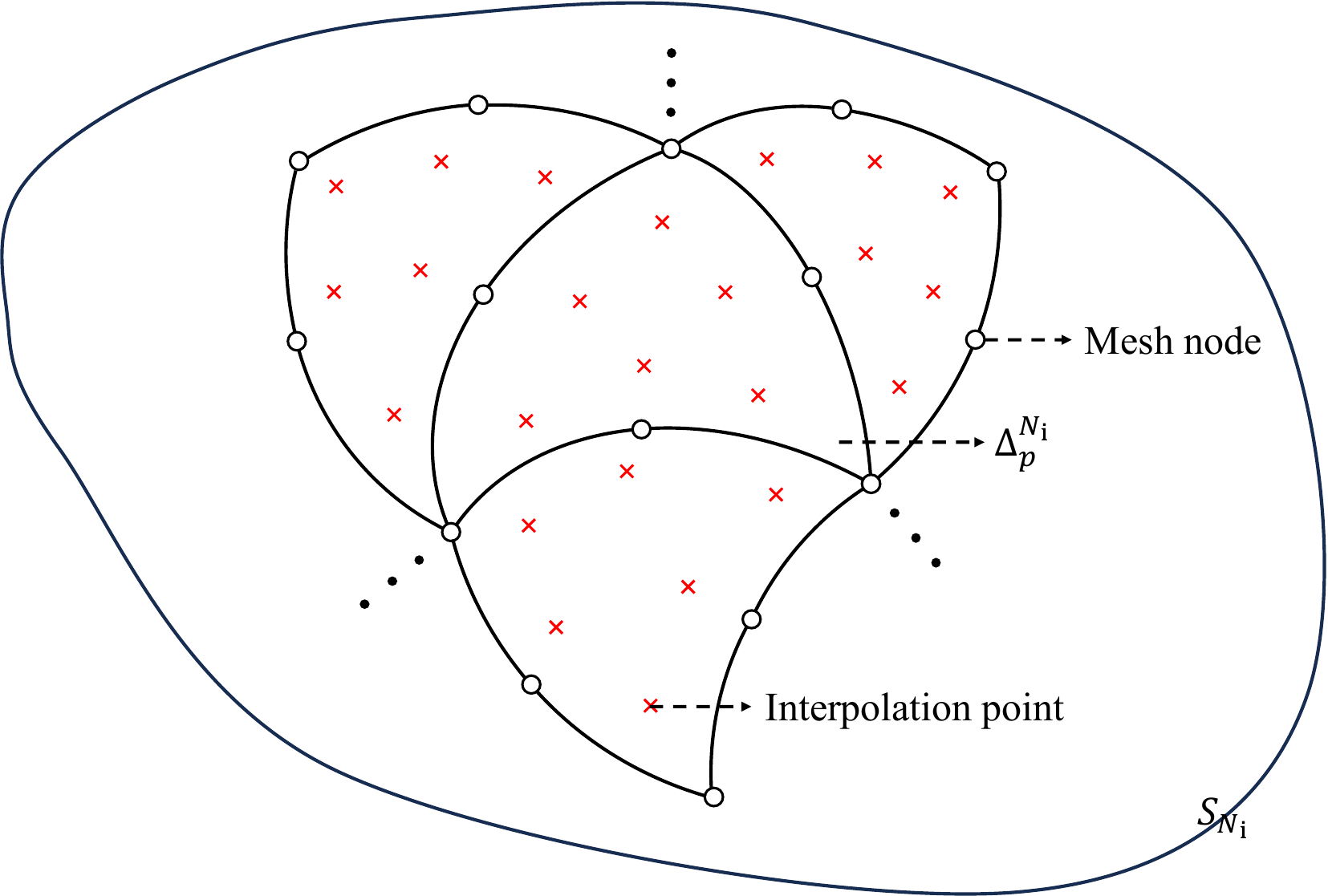}
\caption{An example for the discretization of $S_{N_{\mathrm{i}}}$ using the quadratically curved triangular patches with six mesh nodes along the boundary of each patch and using the $2$nd order interpolation function with six interpolation points on each patch.}
\label{fig_2}
\end{figure}

First, the interfaces of the head volume conductor model $S_1, S_2, ..., S_{N_{\mathrm{i}}}$ are discretized using $N_{\mathrm{p}}^1, N_{\mathrm{p}}^2, .., N_{\mathrm{p}}^{N_{\mathrm{i}}}$ quadratically curved triangular patches, respectively, with six mesh nodes along the boundary of each patch (as seen in Fig. \ref{fig_2}). The average length of all the meshed edges is $h$.  

\begin{figure}[!t]
\centerline{\includegraphics[width=0.9\columnwidth]{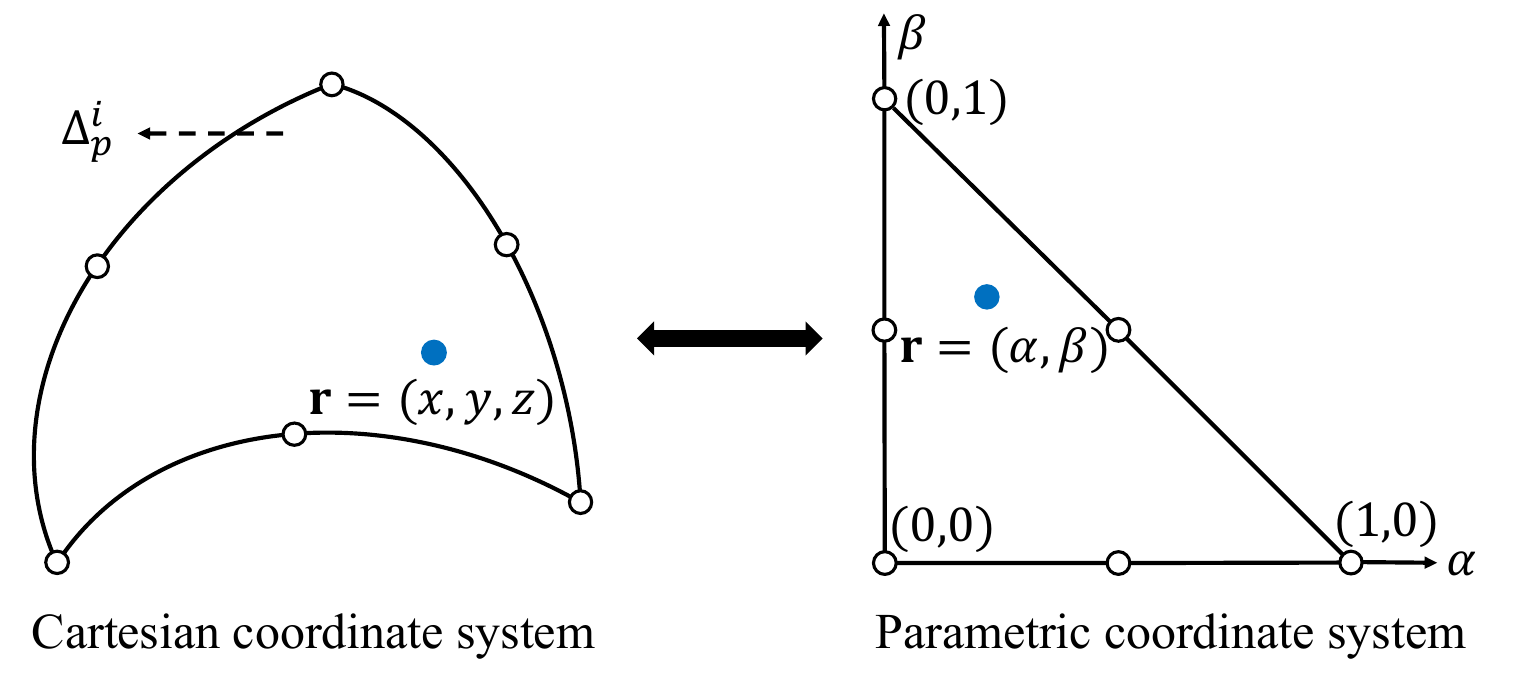}}
\caption{The space mapping between the $(x,y,z)$ Cartesian and  $(\alpha, \beta)$ parametric coordinate systems for the $p^{\mathrm{th}}$ patch on $S_i$.}
\label{fig_3}
\end{figure}

Next, the unknown $\Phi(\mathbf{r})$ of SIEs  on one patch $\Delta_p^i$ can be numerically approximated using the interpolation functions associated with the interpolation points on the patch  as 
\begin{align}
    \Phi(\mathbf{r})=\sum_{a=1}^{N_{\mathrm{a}}} \{ \mathbf{I}_{\Phi} \}_a \vartheta^{-1}(\mathbf{r}) L_a(\mathbf{r}) 
\label{eq_12}
\end{align}
where  $\Phi=\{ V,\xi,V_{\mathrm{ISA}},V_{\mathrm{corr}},J_{N_{\mathrm{i}}-1},J_{N_{\mathrm{i}}} \}$,  $N_{\mathrm{a}}$ is the number of the interpolation points on
each patch (as seen in Fig. \ref{fig_2}), $\vartheta(\mathbf{r})$ represents the Jacobian term at $\mathbf{r}$ resulting from the space mapping between the $(x,y,z)$ Cartesian and $(\alpha, \beta)$ parametric coordinate systems (as seen in Fig. \ref{fig_3}),   $\{\mathbf{I}_{\Phi} \}_a$ is the expansion coefficient associated with the $a$th interpolation point on the patch, and  $L_a(\mathbf{r})$, $a=1,2,...,N_{\mathrm{a}}$ denotes the interpolation functions associated with the  interpolation points on the patch with the property $L_a(\mathbf{r})=1$ when $\mathbf{r}$ is the $a$th interpolation point and $L_a(\mathbf{r})=0$ when $\mathbf{r}$ is the other interpolation point. Moreover, $N_{\mathrm{a}}=\{ 1,3,6 \}$ correspond to the $ \{ 0,1,2 \} $th order interpolation function, respectively \cite{Kang}. In addition, $L_a(\mathbf{r})$, $a=1,2,...,N_{\mathrm{a}}$ can be obtained by solving the Lagrange polynomials-based matrix equation, e.g., when $N_{\mathrm{a}}=6$, which  can be expressed as 
\begin{align}
    \left[
    \begin{array}{cccc}
        1 & 1 & ... & 1   \\
        \alpha_1 & \alpha_2 & ... & \alpha_6   \\
        \beta_1 & \beta_2 & ... & \beta_6   \\
        \alpha_1^2 & \alpha_2^2 & ... & \alpha_6^2  \\
        \alpha_1 \beta_1 & \alpha_2 \beta_2 & ... & \alpha_6 \beta_6   \\
        \beta_1^2 & \beta_2^2 & ... & \beta_6^2   
    \end{array}
    \right]
    \left[
    \begin{array}{c}
        L_1(\mathbf{r})\\
        L_2(\mathbf{r})\\
        L_3(\mathbf{r})\\
        L_4(\mathbf{r})\\
        L_5(\mathbf{r})\\
        L_6(\mathbf{r})
    \end{array}
    \right]=
    \left[
    \begin{array}{c}
        1\\
        \alpha\\
        \beta\\
        \alpha^2\\
        \alpha \beta\\
        \beta^2
    \end{array}
    \right]
\label{eq_13}
\end{align}
where $(\alpha,\beta)$ and $(\alpha_a,\beta_a)$ are coordinates of $\mathbf{r}$ and the $a$th interpolation point on $\Delta_p^i$ in the parametric system, respectively.  

Then, inserting the approximation of the unknowns in \eqref{eq_12} into the SIEs in Section \ref{subsec:SIEs} and point-testing the resulting equations at the observation points (chosen to be the same as the interpolation points)  yields the fully discretized  matrix equation systems that can be solved for the unknown expansion coefficients $\{ \mathbf{I}_{\Phi}  \}_a$. Note that we choose the $N_{\mathrm{a}}$ interpolation points on each patch on the same locations  as the $N_{\mathrm{a}}$ Gaussian quadrature points on the patch, positioned according to the $N_{\mathrm{a}}$-point Gaussian quadrature rule  \cite{Jin} for the Nyström method used in this work. 

In the following section, the discretization of the four SIEs in Section \ref{subsec:SIEs} for analyzing the EEG forward problem  using the Nyström method are presented.

\subsubsection{DL Formulation}
For the DL approach, the unknown $V(\mathbf{r})$ at $\mathbf{r} \in S_1, S_2, ..., S_{N_{\mathrm{i}}}$ is expanded as
\begin{align}
V(\mathbf{r})&=\sum_{i=1}^{N_{\mathrm{i}}}\sum_{p=1}^{N_{\mathrm{p}}^i}\sum_{a=1}^{N_{\mathrm{a}}}\{\mathbf{I}_V \}_{(a,p,i)}\vartheta^{-1}(\mathbf{r})L_{(a,p,i)}(\mathbf{r}).
\label{eq_14}
\end{align}
Inserting \eqref{eq_14} into \eqref{eq_2} and point-testing the resulting equation at the observation points $\mathbf{r}_{(b,q,j)}$, $j=1,2,..,N_{\mathrm{i}}$,  $q=1,2,...,N_{\mathrm{p}}^j$, $b=1,2,...,N_{\mathrm{a}}$ yields the matrix equation system 
\begin{align}
\mathbf{Z}_V\mathbf{I}_V &=\mathbf{V}_V
\label{eq_15}
\end{align}
where the entries of $\mathbf{V}_V$ with the dimension of $N_{\mathrm{s}}^{V} \times 1$ and $\mathbf{Z}_V$ with the  dimension of $N_{\mathrm{s}}^{V} \times N_{\mathrm{s}}^{V}$ are given by
\begin{align}
&\{ \mathbf{V}_V \}_{(b,q,j)}=V^{\mathrm{inc}}(\mathbf{r}_{(b,q,j)})\label{eq_16}\\
\nonumber &\{ \mathbf{Z}_V \}_{(b,q,j)(a,p,i)}=\frac{\sigma_j+\sigma_{j+1}}{2}\vartheta^{-1}(\mathbf{r}_{(b,q,j)})\delta_{ab}\delta_{pq}\delta_{ij} -  \\
& (\sigma_{i+1}-\sigma_i)  \mathrm{P. V.} \int_{\Delta_p^i}\partial_{\mathbf{n}'}G(\mathbf{r}_{(b,q,j)},\mathbf{r}')\vartheta^{-1}(\mathbf{r}')L_{(a,p,i)}(\mathbf{r}')d\mathbf{r}'
\label{eq_17}
\end{align}
where $N_{\mathrm{s}}^{V}=N_{\mathrm{a}} \sum_{i=1}^{N_{\mathrm{i}}}  N_{\mathrm{p}}^i $, $\delta_{ab}=1$ for $a=b$ and $\delta_{ab}=0$ for $a \ne b$, $\Delta_p^i$ denotes the $p$th quadratically curved triangular patch on $S_i$. After the solution of \eqref{eq_15} for $\mathbf{I}_V$, $V(\mathbf{r})$ at $\mathbf{r}$ on one patch of $S_1, ..., S_{N_{\mathrm{i}}}$  can be obtained using \eqref{eq_12}.

\subsubsection{ADL Formulation}
For the ADL approach, the unknown $\xi(\mathbf{r})$ at $\mathbf{r} \in S_1, S_2, ..., S_{N_{\mathrm{i}}}$ is expanded as
\begin{align}
\xi(\mathbf{r})&=\sum_{i=1}^{N_{\mathrm{i}}}\sum_{p=1}^{N_{\mathrm{p}}^i}\sum_{a=1}^{N_{\mathrm{a}}}\{\mathbf{I}_{\xi} \}_{(a,p,i)}\vartheta^{-1}(\mathbf{r})L_{(a,p,i)}(\mathbf{r}).
\label{eq_18}
\end{align}
Inserting \eqref{eq_18} into \eqref{eq_3} and point-testing the resulting equation at the observation points $\mathbf{r}_{(b,q,j)}$, $j=1,2,..,N_{\mathrm{i}}$,  $q=1,2,...,N_{\mathrm{p}}^j$, $b=1,2,...,N_{\mathrm{a}}$ yields the matrix equation system 
\begin{align}
\mathbf{Z}_{\xi}\mathbf{I}_{\xi} &=\mathbf{V}_{\xi} 
\label{eq_19}
\end{align}
where the entries of $\mathbf{V}_{\xi}$ in the dimension of $N_{\mathrm{s}}^{\xi} \times 1$ and $\mathbf{Z}_{\xi}$ in the dimension of $N_{\mathrm{s}}^{\xi} \times N_{\mathrm{s}}^{\xi}$ are given by
\begin{align}
\nonumber &\{\mathbf{V}_{\xi}\}_{(b,q,j)}=\frac{1}{4\pi \sigma_1}\hat{\mathbf{n}}(\mathbf{r}_{(b,q,j)})\\
&\cdot \{\frac{\mathbf{q}}{|\mathbf{r}_{(b,q,j)}-\mathbf{r}_0 |^3}-\frac{3[\mathbf{q}\cdot (\mathbf{r}_{(b,q,j)}-\mathbf{r}_0)](\mathbf{r}_{(b,q,j)}-\mathbf{r}_0)}{|\mathbf{r}_{(b,q,j)}-\mathbf{r}_0 |^5}\}\label{eq_20}\\
\nonumber &\{ \mathbf{Z}_{\xi} \}_{(b,q,j)(a,p,i)}=\frac{\sigma_j+\sigma_{j+1}}{2(\sigma_{j+1}-\sigma_j)}\vartheta^{-1}(\mathbf{r}_{(b,q,j)})\delta_{ab}\delta_{pq}\delta_{ij}\\
&-\mathrm{P. V.}\int_{\Delta_p^i}\partial_{\mathbf{n}}G(\mathbf{r}_{(b,q,j)},\mathbf{r}')\vartheta^{-1}(\mathbf{r}')L_{(a,p,i)}(\mathbf{r}')d\mathbf{r}' 
\label{eq_21}
\end{align}
where $N_{\mathrm{s}}^{\xi}=N_{\mathrm{a}} \sum_{i=1}^{N_{\mathrm{i}}}  N_{\mathrm{p}}^i $. After the solution of \eqref{eq_18} for $\mathbf{I}_{\xi}$, $V(\mathbf{r})$  at the arbitrary point $\mathbf{r}$ can be obtained using \eqref{eq_4} as 
\begin{align}
    \nonumber V(\mathbf{r})=&\frac{V^{\mathrm{inc}}(\mathbf{r})}{\sigma_1}+\sum_{i=1}^{N_{\mathrm{i}}}\sum_{p=1}^{N_{\mathrm{p}}^i}\sum_{a=1}^{N_{\mathrm{a}}}\{\mathbf{I}_{\xi} \}_{(a,p,i)} \times \\
    &\int_{\Delta_p^i}G(\mathbf{r},\mathbf{r}')\vartheta^{-1}(\mathbf{r}')L_{(a,p,i)}(\mathbf{r}')d\mathbf{r}'.\label{eq_22}
\end{align}

\subsubsection{ISADL Formulation}
For ISA for the DL approach, first, the unknown $V_{\mathrm{ISA}}(\mathbf{r})$ at $\mathbf{r} \in S_1, ..., S_{K-1}$ is expanded as
\begin{align}
V_{\mathrm{ISA}}(\mathbf{r})&=\sum_{i=1}^{K-1}\sum_{p=1}^{N_{\mathrm{p}}^i}\sum_{a=1}^{N_{\mathrm{a}}}\{\mathbf{I}_{\mathrm{ISA}} \}_{(a,p,i)}\vartheta^{-1}(\mathbf{r})L_{(a,p,i)}(\mathbf{r}) .
\label{eq_23}
\end{align}
Inserting \eqref{eq_23} into \eqref{eq_6} and point-testing the resulting equation at  $\mathbf{r}_{(b,q,j)}$, $j=1,2,..,K-1$,  $q=1,2,...,N_{\mathrm{p}}^j$, $b=1,2,...,N_{\mathrm{a}}$ yields the following matrix equation system 
\begin{align}
\mathbf{Z}_{\mathrm{ISA}}\mathbf{I}_{\mathrm{ISA}} &=\mathbf{V}_{\mathrm{ISA}} 
\label{eq_24}
\end{align}
where the entries of $\mathbf{V}_{\mathrm{ISA}}$ with the dimension of $N_{\mathrm{s}}^{\mathrm{ISA}} \times 1 $ and $\mathbf{Z}_{\mathrm{ISA}}$ with the dimension of $N_{\mathrm{s}}^{\mathrm{ISA}} \times N_{\mathrm{s}}^{\mathrm{ISA}}$ are given by
\begin{align}
&\{ \mathbf{V}_{\mathrm{ISA}} \}_{(b,q,j)}=V^{\mathrm{inc}}(\mathbf{r}_{(b,q,j)})\label{eq_25}\\
\nonumber &\{ \mathbf{Z}_{\mathrm{ISA}} \}_{(b,q,j)(a,p,i)}=\frac{\sigma_j+\sigma_{j+1}}{2}\vartheta^{-1}(\mathbf{r}_{(b,q,j)})\delta_{ab}\delta_{pq}\delta_{ij}-\\
&(\tilde{\sigma}_{i+1}-\sigma_i)\mathrm{P. V.}\int_{\Delta_p^i}\partial_{\mathbf{n}'}G(\mathbf{r}_{(b,q,j)},\mathbf{r}')\vartheta^{-1}(\mathbf{r}')L_{(a,p,i)}(\mathbf{r}')d\mathbf{r}'
\label{eq_26}
\end{align}
where $\tilde{\sigma}_{i+1}=\sigma_{i+1}$ for $i=1,2,..., K-2$ and $\tilde{\sigma}_{i+1}=0$ for $i=K-1$. 

Then, the unknown $V_{\mathrm{corr}}(\mathbf{r})$ at $\mathbf{r} \in S_1, ..., S_{N_{\mathrm{i}}}$ is written as
\begin{align}
V_{\mathrm{corr}}(\mathbf{r})&=\sum_{i=1}^{N_{\mathrm{i}}}\sum_{p=1}^{N_{\mathrm{p}}^i}\sum_{a=1}^{N_{\mathrm{a}}}\{\mathbf{I}_{\mathrm{corr}} \}_{(a,p,i)}\vartheta^{-1}(\mathbf{r})L_{(a,p,i)}(\mathbf{r}).
\label{eq_27}
\end{align}
Inserting \eqref{eq_27} into \eqref{eq_7} and point-testing the resulting equation at  $\mathbf{r}_{(b,q,j)}$, $j=1,2,..,N_{\mathrm{i}}$,  $q=1,2,...,N_{\mathrm{p}}^j$, $b=1,2,...,N_{\mathrm{a}}$ yields the following matrix equation system 
\begin{align}
\mathbf{Z}_{\mathrm{corr}}\mathbf{I}_{\mathrm{corr}} &=\mathbf{V}_{\mathrm{corr}} 
\label{eq_28}
\end{align}
where the entries of $\mathbf{V}_{\mathrm{corr}}$ with the dimension of $N_{\mathrm{s}}^{\mathrm{corr}} \times 1$ and $\mathbf{Z}_{\mathrm{corr}}$ with the dimension of $N_{\mathrm{s}}^{\mathrm{corr}} \times N_{\mathrm{s}}^{\mathrm{corr}}$ are given by
\begin{align}
\nonumber &\{ \mathbf{V}_{\mathrm{corr}} \}_{(b,q,j)}=\\
&\begin{cases}
    \frac{\sigma_K}{\sigma_{K-1}}V^{\mathrm{inc}}(\mathbf{r}_{(b,q,j)})+
    \sum_{i=1}^{K-2} \frac{\sigma_K (\sigma_{i+1}-\sigma_i)}{\sigma_{K-1}} \times \\  \ \ \int_{S_i}\partial_{\mathbf{n}'}G(\mathbf{r}_{(b,q,j)},\mathbf{r}')V_{\mathrm{ISA}}(\mathbf{r}')d\mathbf{r}' , \  \ \mathrm{if} \  j=K, K+1, ..., N_{\mathrm{i}} \\
    \frac{\sigma_K}{\sigma_{K-1}}V^{\mathrm{inc}}(\mathbf{r}_{(b,q,j)})-\sigma_K V_{\mathrm{ISA}}(\mathbf{r}_{(b,q,j)})+\sum_{i=1}^{K-2} \frac{\sigma_K (\sigma_{i+1}-\sigma_i)}{\sigma_{K-1}}\\
        \ \  \times \int_{S_i}\partial_{\mathbf{n}'}G(\mathbf{r}_{(b,q,j)},\mathbf{r}')V_{\mathrm{ISA}}(\mathbf{r}')d\mathbf{r}' ,  \ \ \mathrm{if} \ j=K-1 \\
    \frac{\sigma_K}{\sigma_{K-1}}V^{\mathrm{inc}}(\mathbf{r}_{(b,q,j)})-\frac{\sigma_K (\sigma_j+\sigma_{j+1})}{2 \sigma_{K-1}} V_{\mathrm{ISA}}(\mathbf{r}_{(b,q,j)}) +\sum_{i=1}^{K-2} \\
    \ \ \frac{\sigma_K (\sigma_{i+1}-\sigma_i)}{\sigma_{K-1}}   \mathrm{P. V.} \int_{S_i}\partial_{\mathbf{n}'}G(\mathbf{r}_{(b,q,j)},\mathbf{r}')V_{\mathrm{ISA}}(\mathbf{r}')d\mathbf{r}' , \\
   \ \  \mathrm{if} \ j=1, ..., K-2 
\end{cases} \label{eq_29} \\
\nonumber &\{ \mathbf{Z}_{\mathrm{corr}} \}_{(b,q,j)(a,p,i)}=\frac{\sigma_j+\sigma_{j+1}}{2}\vartheta^{-1}(\mathbf{r}_{(b,q,j)})\delta_{ab}\delta_{pq}\delta_{ij}-\\
& (\sigma_{i+1}-\sigma_i) \mathrm{P. V.} \int_{\Delta_p^i}\partial_{\mathbf{n}'}G(\mathbf{r}_{(b,q,j)},\mathbf{r}')\vartheta^{-1}(\mathbf{r}')L_{(a,p,i)}(\mathbf{r}')d\mathbf{r}'.\label{eq_30}
\end{align}

After the solution of \eqref{eq_24} and \eqref{eq_28} for $\mathbf{I}_{\mathrm{ISA}}$ and $\mathbf{I}_{\mathrm{corr}}$, respectively, $V(\mathbf{r})$  at $\mathbf{r}$ on one patch of the interfaces $S_1, S_2, ..., S_{N_{\mathrm{i}}}$ can be obtained using \eqref{eq_5} as 
\begin{align}
V (\mathbf{r})=
\begin{cases}
\sum_{a=1}^{N_{\mathrm{a}}}\{\{\mathbf{I}_{\mathrm{ISA}} \}_a +\{\mathbf{I}_{\mathrm{corr}} \}_a \}\vartheta^{-1}(\mathbf{r})L_a(\mathbf{r}), \\ 
\ \ \mathrm{if} \ \mathbf{r} \in S_1,...,S_{K-1}\\
\sum_{a=1}^{N_{\mathrm{a}}}\{\mathbf{I}_{\mathrm{corr}} \}_a\vartheta^{-1}(\mathbf{r})L_a(\mathbf{r}), 
\ \ \mathrm{if} \ \mathbf{r} \in S_K,...,S_{N_{\mathrm{i}}}
\end{cases}
\label{eq_31}
\end{align}

\subsubsection{IADL Formulation}
For the IADL approach, the unknown $J_{N_{\mathrm{i}}-1}(\mathbf{r})$ and $J_{N_{\mathrm{i}}}(\mathbf{r})$  at $\mathbf{r} \in S_{N_{\mathrm{i}}-1}$ and $\mathbf{r} \in S_{N_{\mathrm{i}}}$ are  expanded respectively as
\begin{align}
J_{N_{\mathrm{i}}-1}(\mathbf{r})&=\sum_{p=1}^{N_{\mathrm{p}}^{N_{\mathrm{i}}-1}}\sum_{a=1}^{N_{\mathrm{a}}}\{\mathbf{I}_{N_{\mathrm{i}}-1} \}_{(a,p,N_{\mathrm{i}}-1)}\vartheta^{-1}(\mathbf{r})L_{(a,p,N_{\mathrm{i}}-1)}(\mathbf{r}) \label{eq_32} \\
J_{N_{\mathrm{i}}}(\mathbf{r})&=\sum_{p=1}^{N_{\mathrm{p}}^{N_{\mathrm{i}}}}\sum_{a=1}^{N_{\mathrm{a}}}\{\mathbf{I}_{N_{\mathrm{i}}} \}_{(a,p,N_{\mathrm{i}})}\vartheta^{-1}(\mathbf{r})L_{(a,p,N_{\mathrm{i}})}(\mathbf{r}). \label{eq_33}
\end{align}
Inserting \eqref{eq_32}-\eqref{eq_33} into \eqref{eq_10}-\eqref{eq_11} and point-testing the resulting equations at  $\mathbf{r}_{(b,q,N_{\mathrm{i}}-1)}$,   $q=1,2,...,N_{\mathrm{p}}^{N_{\mathrm{i}}-1}$, $b=1,2,...,N_{\mathrm{a}}$ and $\mathbf{r}_{(b,q,N_{\mathrm{i}})}$,   $q=1,2,...,N_{\mathrm{p}}^{N_{\mathrm{i}}}$, $b=1,2,...,N_{\mathrm{a}}$, respectively, yields the following matrix equation system 
\begin{align}
\left[
\begin{array}{cc}
    \mathbf{Z}_{N_{\mathrm{i}}-1,N_{\mathrm{i}}-1} & \mathbf{Z}_{N_{\mathrm{i}}-1,N_{\mathrm{i}}} \\
    \mathbf{Z}_{N_{\mathrm{i}},N_{\mathrm{i}}-1} & \mathbf{Z}_{N_{\mathrm{i}},N_{\mathrm{i}}}
\end{array}
\right]
\left[
\begin{array}{c}
    \mathbf{I}_{N_{\mathrm{i}}-1} \\
    \mathbf{I}_{N_{\mathrm{i}}}
\end{array}
\right]=
\left[
\begin{array}{c}
    \mathbf{V}_{N_{\mathrm{i}}-1} \\
    \mathbf{V}_{N_{\mathrm{i}}}
\end{array}
\right]
\label{eq_34}
\end{align}
where the entries of $\mathbf{V}_{N_{\mathrm{i}}-1}$ with the dimension of $N_{\mathrm{s}}^{N_{\mathrm{i}}-1} \times 1$, $\mathbf{V}_{N_{\mathrm{i}}}$ with the dimension of $N_{\mathrm{s}}^{N_{\mathrm{i}}-1} \times 1$, $\mathbf{Z}_{N_{\mathrm{i}}-1,N_{\mathrm{i}}-1}$ with the dimension of $N_{\mathrm{s}}^{N_{\mathrm{i}}-1} \times N_{\mathrm{s}}^{N_{\mathrm{i}}-1}$, $\mathbf{Z}_{N_{\mathrm{i}}-1,N_{\mathrm{i}}}$ with the dimension of $N_{\mathrm{s}}^{N_{\mathrm{i}}-1} \times N_{\mathrm{s}}^{N_{\mathrm{i}}}$, $\mathbf{Z}_{N_{\mathrm{i}},N_{\mathrm{i}}-1}$ with the dimension of $N_{\mathrm{s}}^{N_{\mathrm{i}}} \times N_{\mathrm{s}}^{N_{\mathrm{i}}-1}$, and $\mathbf{Z}_{N_{\mathrm{i}},N_{\mathrm{i}}}$ with the dimension of $N_{\mathrm{s}}^{N_{\mathrm{i}}} \times N_{\mathrm{s}}^{N_{\mathrm{i}}}$ are given by
\begin{align}
&\{ \mathbf{V}_{N_{\mathrm{i}}-1} \}_{(b,q,{N_{\mathrm{i}}-1})}=\frac{\sigma_{N_{\mathrm{i}}-1}}{\sigma_{N_{\mathrm{i}}}-\sigma_{N_{\mathrm{i}}-1}}\xi(\mathbf{r}_{(b,q,N_{\mathrm{i}}-1)}) \label{eq_35} \\
&\{ \mathbf{V}_{N_{\mathrm{i}}} \}_{(b,q,{N_{\mathrm{i}}})}=0 \label{eq_36} \\
\nonumber &\{ \mathbf{Z}_{N_{\mathrm{i}}-1,N_{\mathrm{i}}-1} \}_{(b,q,N_{\mathrm{i}}-1)(a,p,N_{\mathrm{i}}-1)}=-\frac{1}{2}\vartheta^{-1}(\mathbf{r}_{(b,q,N_{\mathrm{i}}-1)})\delta_{ab}\delta_{pq}\\
&\ +\mathrm{P. V.}\int_{\Delta_p^{N_{\mathrm{i}}-1}}\partial_{\mathbf{n}}G(\mathbf{r}_{(b,q,N_{\mathrm{i}}-1)},\mathbf{r}')\vartheta^{-1}(\mathbf{r}')L_{(a,p,N_{\mathrm{i}}-1)}(\mathbf{r}')d\mathbf{r}' \label{eq_37} \\
\nonumber &\{ \mathbf{Z}_{N_{\mathrm{i}}-1,N_{\mathrm{i}}} \}_{(b,q,N_{\mathrm{i}}-1)(a,p,N_{\mathrm{i}})}=\\
&  \ \ \ \int_{\Delta_p^{N_{\mathrm{i}}}}\partial_{\mathbf{n}}G(\mathbf{r}_{(b,q,N_{\mathrm{i}}-1)},\mathbf{r}')\vartheta^{-1}(\mathbf{r}')L_{(a,p,N_{\mathrm{i}})}(\mathbf{r}')d\mathbf{r}' \label{eq_38} \\
\nonumber &\{ \mathbf{Z}_{N_{\mathrm{i}},N_{\mathrm{i}}-1} \}_{(b,q,N_{\mathrm{i}})(a,p,N_{\mathrm{i}}-1)}=\\
&  \ \ \ \int_{\Delta_p^{N_{\mathrm{i}}-1}}\partial_{\mathbf{n}}G(\mathbf{r}_{(b,q,N_{\mathrm{i}})},\mathbf{r}')\vartheta^{-1}(\mathbf{r}')L_{(a,p,N_{\mathrm{i}}-1)}(\mathbf{r}')d\mathbf{r}' \label{eq_39} \\
\nonumber &\{ \mathbf{Z}_{N_{\mathrm{i}},N_{\mathrm{i}}} \}_{(b,q,N_{\mathrm{i}})(a,p,N_{\mathrm{i}})}=\frac{1}{2}\vartheta^{-1}(\mathbf{r}_{(b,q,N_{\mathrm{i}})})\delta_{ab}\delta_{pq}+\\
& \ \ \ \mathrm{P. V.}\int_{\Delta_p^{N_{\mathrm{i}}}}\partial_{\mathbf{n}}G(\mathbf{r}_{(b,q,N_{\mathrm{i}})},\mathbf{r}')\vartheta^{-1}(\mathbf{r}')L_{(a,p,N_{\mathrm{i}})}(\mathbf{r}')d\mathbf{r}' .\label{eq_40}
\end{align}
After the solution of \eqref{eq_34} for $\mathbf{I}_{N_{\mathrm{i}}-1}$ and $\mathbf{I}_{N_{\mathrm{i}}}$, $V(\mathbf{r})$  at $\mathbf{r} \in \Omega_{N_{\mathrm{i}}}$ can be obtained using \eqref{eq_9} as 
\begin{align}
\nonumber &V(\mathbf{r})=\sum_{p=1}^{N_{\mathrm{p}}^{N_{\mathrm{i}}-1}}\sum_{a=1}^{N_{\mathrm{a}}}\{\mathbf{I}_{N_{\mathrm{i}}-1} \}_{(a,p,N_{\mathrm{i}}-1)}\times \\
\nonumber &\int_{\Delta_p^{N_{\mathrm{i}}-1}} G(\mathbf{r},\mathbf{r}') \vartheta^{-1}(\mathbf{r}') L_{(a,p,N_{\mathrm{i}}-1)}(\mathbf{r}') d
 \mathbf{r}' + \sum_{p=1}^{N_{\mathrm{p}}^{N_{\mathrm{i}}}}\sum_{a=1}^{N_{\mathrm{a}}} \\
 &\{\mathbf{I}_{N_{\mathrm{i}}} \}_{(a,p,N_{\mathrm{i}})} 
  \int_{\Delta_p^{N_{\mathrm{i}}}} G(\mathbf{r},\mathbf{r}') \vartheta^{-1}(\mathbf{r}') L_{(a,p,N_{\mathrm{i}})}(\mathbf{r}') d
 \mathbf{r}', \  \mathbf{r} \in \Omega_{N_{\mathrm{i}}} .\label{eq_41}
\end{align}

\subsection{Numerical Evaluation of the Surface Integrals}
\label{subsec:evaluation}
As seen in \eqref{eq_17}, \eqref{eq_21}-\eqref{eq_22}, \eqref{eq_26}, \eqref{eq_29}-\eqref{eq_30}, and \eqref{eq_37}-\eqref{eq_41}, there are three kinds of surface integrals as given by 
\begin{align}
&  I_1=\int_{\Delta_p^i}\partial_{\mathbf{n}'}G(\mathbf{r}_{(b,q,j)},\mathbf{r}')\vartheta^{-1}(\mathbf{r}')L_{(a,p,i)}(\mathbf{r}')d\mathbf{r}' \label{eq_42} \\
&I_2=\int_{\Delta_p^i}\partial_{\mathbf{n}}G(\mathbf{r}_{(b,q,j)},\mathbf{r}')\vartheta^{-1}(\mathbf{r}')L_{(a,p,i)}(\mathbf{r}')d\mathbf{r}' \label{eq_43}  \\
&I_3=\int_{\Delta_p^i}G(\mathbf{r},\mathbf{r}')\vartheta^{-1}(\mathbf{r}')L_{(a,p,i)}(\mathbf{r}')d\mathbf{r}'. \label{eq_44}
\end{align}
Depending on the distance $d$ between the observation point ($\mathbf{r}_{(b,q,j)}$ or $\mathbf{r}$) and the source triangular patch $\Delta_p^i$, three following cases are considered for the evaluation of the surface integrals $I_1$, $I_2$, and $I_3$
in \eqref{eq_42}-\eqref{eq_44}.

\subsubsection{Far-interaction}

When the observation point ($\mathbf{r}_{(b,q,j)}$ or $\mathbf{r}$) is not located on the source triangular patch $\Delta_p^i$ and far away from
 $\Delta_p^i$ with the distance $d>\chi$, where $\chi$ is a threshold value, the surface integrals $I_1$, $I_2$, and $I_3$ in \eqref{eq_42}-\eqref{eq_44} can be calculated by first mapping the integrals in the Cartesian coordinate system to those in the parametric coordinate system and then using the $N_{\mathrm{a}}$-point Gaussian quadrature rule \cite{Jin} for the numerical evaluation of them resulting in the expressions 
\begin{align}
\nonumber I_1 &=\int_{\Delta_p^i}\partial_{\mathbf{n}'}G(\mathbf{r}_{(b,q,j)},\mathbf{r}')\vartheta^{-1}(\mathbf{r}')L_{(a,p,i)}(\mathbf{r}')d\mathbf{r}' \\
    &= \omega_{(a,p,i)}  \hat{\mathbf{n}}'(\mathbf{r}'_{(a,p,i)}) \cdot \frac{\mathbf{r}_{(b,q,j)}-\mathbf{r}'_{(a,p,i)}}{8\pi |\mathbf{r}_{(b,q,j)}-\mathbf{r}'_{(a,p,i)}|^3} \label{eq_45}\\
\nonumber I_2&=\int_{\Delta_p^i}\partial_{\mathbf{n}}G(\mathbf{r}_{(b,q,j)},\mathbf{r}')\vartheta^{-1}(\mathbf{r}')L_{(a,p,i)}(\mathbf{r}')d\mathbf{r}' \\
    &= \omega_{(a,p,i)}   \hat{\mathbf{n}}(\mathbf{r}_{(b,q,j)}) \cdot \frac{\mathbf{r}'_{(a,p,i)}-\mathbf{r}_{(b,q,j)}}{8\pi |\mathbf{r}_{(b,q,j)}-\mathbf{r}'_{(a,p,i)}|^3} \label{eq_46}\\
\nonumber I_3&=\int_{\Delta_p^i}G(\mathbf{r},\mathbf{r}')\vartheta^{-1}(\mathbf{r}')L_{(a,p,i)}(\mathbf{r}')d\mathbf{r}' \\
    &= \omega_{(a,p,i)}    \frac{1}{8\pi |\mathbf{r}-\mathbf{r}'_{(a,p,i)}|} 
\label{eq_47}
\end{align}
where $\mathbf{r}'_{(a,p,i)}$ is the source Gaussian quadrature point same as the $a$th interpolation point on the $p^{\mathrm{th}}$ patch of the $i^{\mathrm{th}}$ interface, and $\omega_{(a,p,i)}$ is the weight coefficient associated with  $\mathbf{r}'_{(a,p,i)}$. 

Clearly, the computation of $I_1$, $I_2$, and $I_3$  is  simplified for the far-interaction case  due to the flexibility of the choice of the interpolation points using the proposed scheme.

\subsubsection{Near-interaction}
When   the observation point ($\mathbf{r}_{(b,q,j)}$ or $\mathbf{r}$) is not located on $\Delta_p^i$ and  close to
 $\Delta_p^i$ with the distance $d<\chi$,  after the space mapping from  the Cartesian coordinate system to the parametric coordinate system for the surface integrals $I_1$, $I_2$, and $I_3$ in \eqref{eq_42}-\eqref{eq_44}, we can use the $N_{\mathrm{k}}$-point Gaussian quadrature rule \cite{Jin} for the numerical evaluation of them resulting in the following expressions
\begin{align}
\nonumber I_1&=\int_{\Delta_p^i}\partial_{\mathbf{n}'}G(\mathbf{r}_{(b,q,j)},\mathbf{r}')\vartheta^{-1}(\mathbf{r}')L_{(a,p,i)}(\mathbf{r}')d\mathbf{r}' \\
    &=\sum_{k=1}^{N_{\mathrm{k}}} \omega_k  L_{(a,p,i)}(\mathbf{r}'_k) \hat{\mathbf{n}}'(\mathbf{r}'_k) \cdot \frac{\mathbf{r}_{(b,q,j)}-\mathbf{r}'_k}{8\pi |\mathbf{r}_{(b,q,j)}-\mathbf{r}'_k|^3} \label{eq_48} \\
\nonumber I_2&=\int_{\Delta_p^i}\partial_{\mathbf{n}}G(\mathbf{r}_{(b,q,j)},\mathbf{r}')\vartheta^{-1}(\mathbf{r}')L_{(a,p,i)}(\mathbf{r}')d\mathbf{r}' \\
    &=\sum_{k=1}^{N_{\mathrm{k}}} \omega_k  L_{(a,p,i)}(\mathbf{r}'_k) \hat{\mathbf{n}}(\mathbf{r}_{(b,q,j)}) \cdot \frac{\mathbf{r}'_k-\mathbf{r}_{(b,q,j)}}{8\pi |\mathbf{r}_{(b,q,j)}-\mathbf{r}'_k|^3} \label{eq_49}\\
\nonumber I_3&=\int_{\Delta_p^i}G(\mathbf{r},\mathbf{r}')\vartheta^{-1}(\mathbf{r}')L_{(a,p,i)}(\mathbf{r}')d\mathbf{r}' \\
    &=\sum_{k=1}^{N_{\mathrm{k}}} \omega_k  L_{(a,p,i)}(\mathbf{r}'_k)  \frac{1}{8\pi |\mathbf{r}-\mathbf{r}'_k|} \label{eq_50}
\end{align}
where $N_{\mathrm{k}}>N_{\mathrm{a}}$, $\mathbf{r}'_k$ is the $k^{\mathrm{th}}$ source Gaussian quadrature point, and $\omega_k$ is the weight coefficient associated with  $\mathbf{r}'_k$.

\subsubsection{Self-interaction}
When the observation point ($\mathbf{r}_{(b,q,j)}$ or $\mathbf{r}$) is located on $\Delta_p^i$, the integrands of the surface integrals $I_1$, $I_2$, and $I_3$ in \eqref{eq_42}-\eqref{eq_44} become singular. For this case, direct numerical evaluation of them can not be used anymore due to the large numerical error. Therefore, techniques for the singularity treatment should be applied. Specifically, the integrands of   $I_2$ in  \eqref{eq_43} and $I_3$ in \eqref{eq_44} become weakly-singular while that of $I_1$ in \eqref{eq_42} becomes strongly-singular for the self interaction case.  The weakly-singular surface integrals $I_2$ in  \eqref{eq_43} and $I_3$ in \eqref{eq_44} can be evaluated using the Duffy method \cite{Duffy} while the strongly-singular surface integral $I_1$ in  \eqref{eq_42} can be evaluated using the polar coordinate transformation method \cite{Guiggiani, Chen3}. It should be noted that, the implementations of both the Duffy and polar coordinate transformation methods still require the numerical evaluation of line integrals, which can be addressed using the $N_{\mathrm{v}}$-point Gauss-Legendre quadrature rule \cite{Jin}. For the details of the implementations of the weakly- and strongly-singularity treatment techniques, interested readers can be referred to \cite{Duffy, Guiggiani, Chen3}.

\section{Numerical Results}
\label{sec:results}
In this section, numerical experiments  are presented to demonstrate the  accuracy, flexibility, and efficiency of the proposed scheme for SIEs  for analyzing the EEG forward problem. 

A spherical head volume conductor model is considered in this section since the  electric potential on the exterior surface of the spherical head model can be analytically obtained as the reference solution \cite{Zhang}. We consider a spherical head model  composed of three layers (i.e., $N_{\mathrm{i}}=3$) representing the brain, skull, and scalp tissues from the inner to the outer, respectively. The radii of the spherical interfaces $S_1$, $S_2$, and $S_3$ are $0.087 $, $0.092 $, and $0.1 $ m, respectively. The dipolar moment is $\mathbf{q}=(1,0,1)/\sqrt{2} \ \mathrm{A \cdot m}$. 

To quantify the accuracy of the solution of SIEs, the relative error of the solution is defined by comparing the electric potential on $S_3$ obtained after using the proposed scheme with the analytical reference solution  \cite{Zhang} as
\begin{align}
\mathrm{Relative \  error}=\sqrt{\frac{\sum_{p=1}^{N_{\mathrm{p}}^3}|\tilde{V}^{\mathrm{ref}}(\mathbf{r}_p^{\mathrm{c}})-\tilde{V}^{\mathrm{num}}(\mathbf{r}_p^{\mathrm{c}})  |^2}{\sum_{p=1}^{N_{\mathrm{p}}^3} |\tilde{V}^{\mathrm{ref}}(\mathbf{r}_p^{\mathrm{c}})|^2}}
\label{eq_51}
\end{align}
where $\mathbf{r}_p^{\mathrm{c}}$ denotes the center point of the $p$th quadratically curved triangular patch on $S_3$, ``$\mathrm{ref}$'' denotes the analytical reference solution \cite{Zhang}, ``$\mathrm{num}$'' denotes the numerical solution obtained after using the proposed scheme, and $\tilde{V}^{\mathrm{type}}(\mathbf{r}_p^{\mathrm{c}})$, $\mathrm{type} \in \{\mathrm{ref}, \mathrm{num}\}$ is defined  as
\begin{align}
\tilde{V}^{\mathrm{type}}(\mathbf{r}_p^{\mathrm{c}})=V^{\mathrm{type}}(\mathbf{r}_p^{\mathrm{c}})-\frac{\sum_{q=1}^{N_{\mathrm{p}}^3}V^{\mathrm{type}}(\mathbf{r}_q^{\mathrm{c}})}{N_{\mathrm{p}}^3}.
\label{eq_52}
\end{align}
Note that, using $\tilde{V}^{\mathrm{type}}(\mathbf{r}_p^{\mathrm{c}})$ for the computation of the relative error in \eqref{eq_51} ensures that both the analytical and numerical solutions have zero means prior to the comparison  \cite{Kybic}.

Unless stated otherwise in this section, for the numerical evaluation of the surface integrals in Section \ref{subsec:evaluation}, the distance threshold value is set to be $\chi=0.04 \ \mathrm{m}$, the 16-point Gaussian quadrature rule is used for the surface integrals for the near interaction (i.e., $N_{\mathrm{k}}=16$), and the 15-point Gauss-Legendre quadrature rule is used for the line integrals for the self interaction (i.e., $N_{\mathrm{v}}=15$). A LU decomposition-based direct solver  \cite{Anderson} is used for solving the matrix equation systems  in Section \ref{subsec:discretization}.


As the first numerical experiment, the effect of the mesh density on the accuracy of the  solution of SIEs  using the proposed scheme is investigated. 
The  dipolar source is located at $\mathbf{r}_0=(0.0425, 0, 0) \ \mathrm{m}$. The conductivities of three layers are $\sigma_1=1 \ \mathrm{S/m}$, $\sigma_2=0.025 \ \mathrm{S/m}$, and $\sigma_3=1 \ \mathrm{S/m}$, respectively. The $2$th order interpolation function is used for the Nyström method (i.e.,  $N_{\mathrm{a}}=6$). Six sets of the mesh density are considered as $h= \{0.01, 0.015, 0.02, 0.025, 0.03, 0.035 \} \ \mathrm{m}$ resulting in $N_{\mathrm{p}}^1= \{1918, 854, 456, 300, 208, 136 \} $,  $N_{\mathrm{p}}^2=\{2170, 932, 532, 326, 214, 172 \}$, and $ N_{\mathrm{p}}^3=\{2556, 1120, 626, 400, 282, 202 \}$, respectively. 

\begin{figure}[!t]
\centering
\includegraphics[width=0.9\columnwidth]{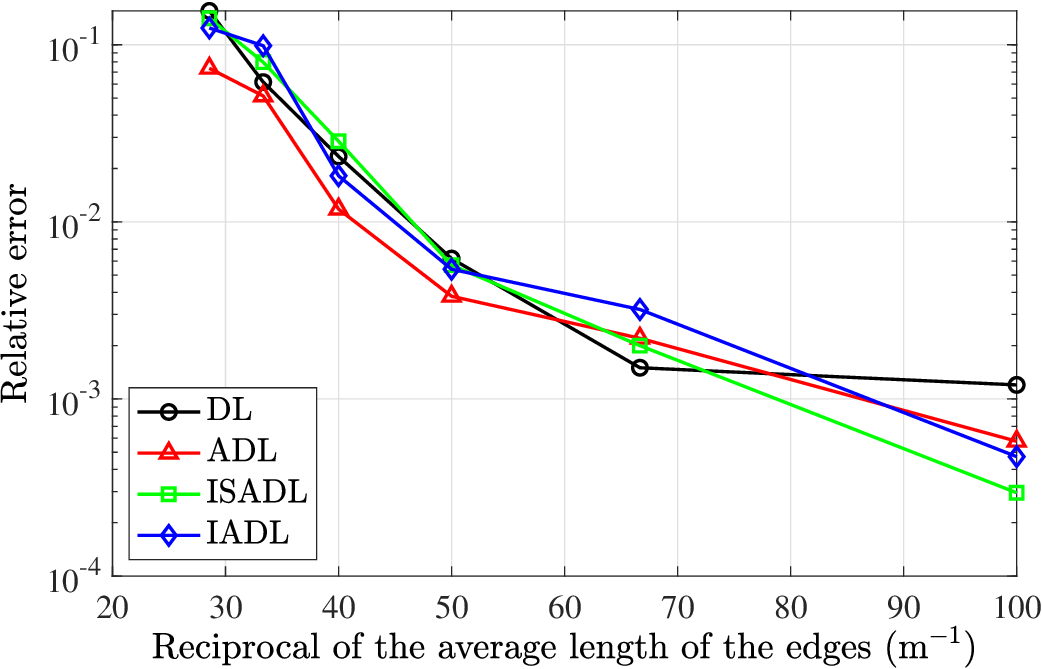}
\caption{ The relative error of the solution after using the proposed scheme versus  the reciprocal of the average length of the edges for the DL, ADL, ISADL, and IADL approaches.}
\label{fig_4}
\end{figure}

Fig. \ref{fig_4} plots the relative error of the solution after using the proposed scheme versus the reciprocal of the average length of the edges for the DL, ADL, ISADL, and IADL approaches. As seen from this figure, the relative error of the solution after using the proposed scheme decreases with the increase of the mesh density for all the four approaches, which demonstrate the efficiency of the proposed scheme. In addition, the ISADL approach has the smallest relative error of the solution for all the six sets of the mesh density among all the four approaches.

For the second numerical experiment, we investigate the effect of the order of the interpolation function on the accuracy of the solution of SIEs  using the proposed scheme.
The  dipolar source is located at $\mathbf{r}_0=(0.0425, 0, 0) \ \mathrm{m}$. The conductivities of three layers are $\sigma_1=1 \ \mathrm{S/m}$, $\sigma_2=0.025 \ \mathrm{S/m}$, and $\sigma_3=1 \ \mathrm{S/m}$, respectively. 
The average edge length of the mesh is   $h=0.01 \ \mathrm{m}$ resulting in $N_{\mathrm{p}}^1=1918$, $N_{\mathrm{p}}^2=2170$, and $N_{\mathrm{p}}^3=2556$. Three sets of the interpolation function are considered as $ \mathrm{order}=\{ 0,  1,   2 \} $th for the Nyström method, i.e.,  $N_{\mathrm{a}}=\{1,  3 ,  6  \}$, respectively.

\begin{figure}[!t]
\centering
\includegraphics[width=0.9\columnwidth]{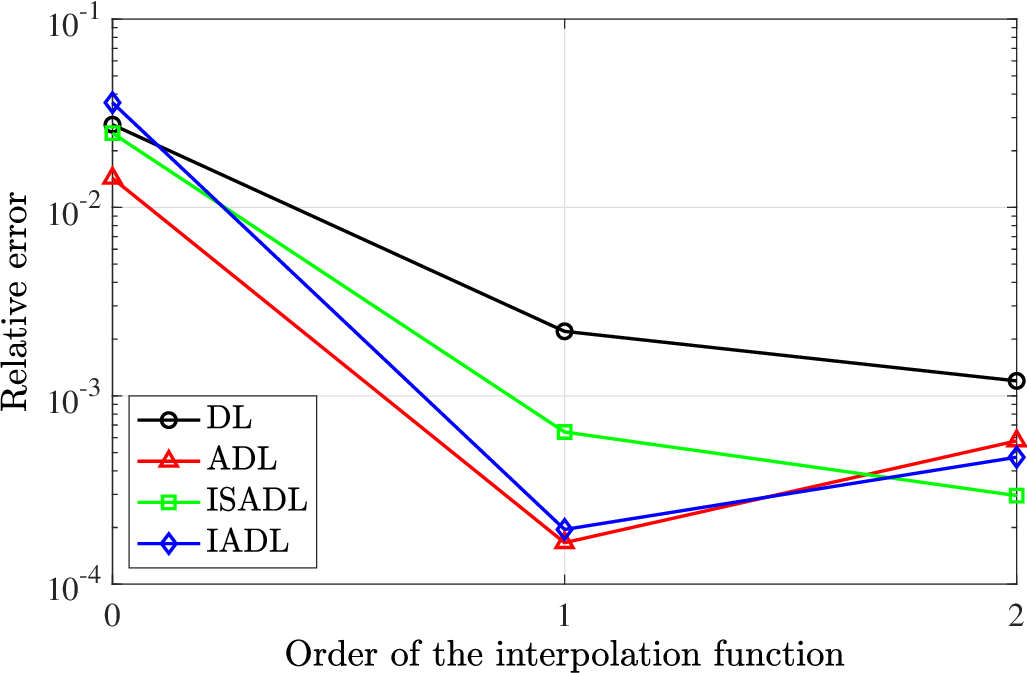}
\caption{ The relative error of the solution after using the proposed scheme  versus the order of the interpolation function for the DL, ADL, ISADL, and IADL approaches.}
\label{fig_5}
\end{figure}

Fig. \ref{fig_5} compares the relative error of the solution after using the proposed scheme  versus the order of the interpolation function for the DL, ADL, ISADL, and IADL approaches.    This figure shows that the accuracy of the  solution of SIEs using the proposed scheme for all the four approaches is satisfactory for most applications. Additionally, the relative error of the solution after using
the proposed scheme decreases with the increase of the order of the interpolation function for all the approaches except for the ADL and IADL approaches. These observations demonstrate the flexibility of the proposed scheme while the saturation behavior of the plots for the ADL and IADL approaches  is due to the fact that the accuracy of the mesh used for this experiment  is limited for these two approaches for the $2$nd order interpolation function. Note that similar saturation behaviors appear in other references \cite{Kang,Chen2} using the Nyström-based HO discretization methods. The usage of higher-order surface patches (e.g., cubic) could be helpful to alleviate the saturation issue, which is out of the scope of this paper.

Next, the effect of the dipole position on the accuracy of the solution of SIEs  using the proposed scheme is investigated.
The $2$th order interpolation function is used for the Nyström method (i.e.,  $N_{\mathrm{a}}=6$). 
The mesh density is $h=0.02 \ \mathrm{m}$ resulting in $N_{\mathrm{p}}^1=456$, $N_{\mathrm{p}}^2=532$, and $N_{\mathrm{p}}^3=626$. Twelve sets of the dipole position located along the $x$-axis are considered as $\mathbf{r}_0=(x_0,0,0)\ \mathrm{m}$ with $x_0=\{0.1, 1, 2, 3, 4.25, 5, 6, 7, 7.5, 8, 8.25, 8.4 \} \times 10^{-2}$. 

\begin{figure}[!t]
\centering
\includegraphics[width=0.9\columnwidth]{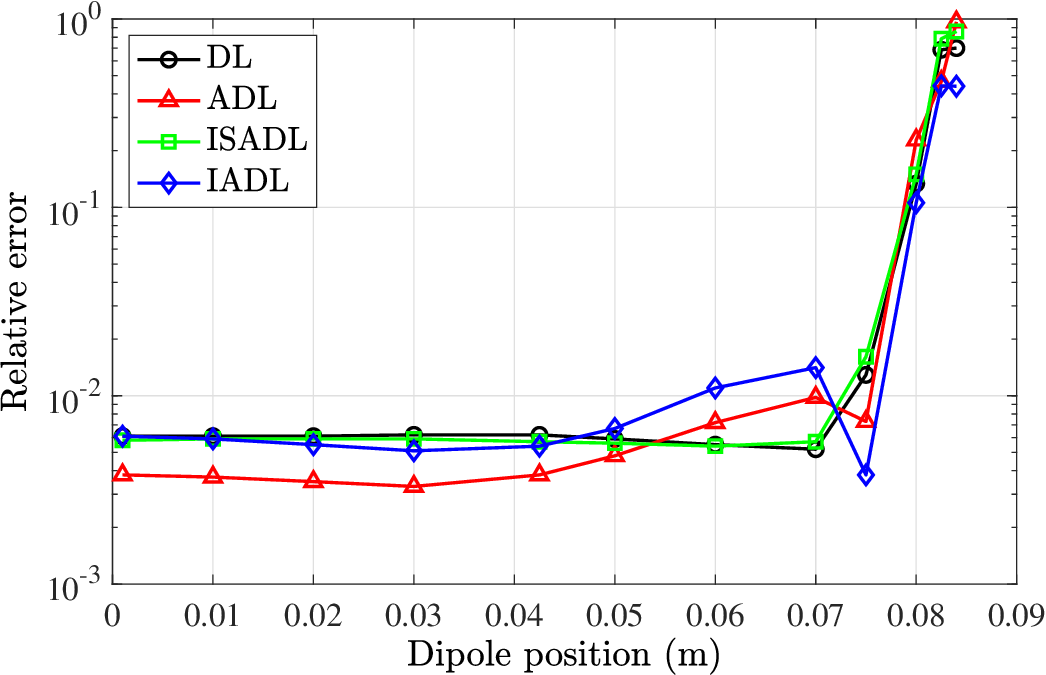}
\caption{ The relative error of the solution after using the proposed scheme  versus the dipole position for the DL, ADL, ISADL, and IADL approaches.}
\label{fig_6}
\end{figure}

Fig. \ref{fig_6} shows the relative error of the solution after using the proposed scheme  versus the dipole position for the DL, ADL, ISADL, and IADL approaches. As we can see from the figure,  the accuracy of the solution of SIEs using the proposed scheme for all the four approaches is  good when $  x_0 \leq 0.075$ while it becomes worse as the dipolar source approaches the boundary of the interior layer $S_1$. This behavior observed from Fig. \ref{fig_6} matches with that observed from Fig. 4 in \cite{Kybic} and Fig. 6 in \cite{Rahmouni}.

As the last numerical experiment, the effect of the  conductivity  of the skull layer on the accuracy of the solution of SIEs using the proposed scheme is investigated. 
The  dipolar source is located at $\mathbf{r}_0=(0.0425, 0, 0) \ \mathrm{m}$. The $2$nd order interpolation function is used for the Nyström method (i.e.,  $N_{\mathrm{a}}=6$). 
The mesh density is $h=0.02 \ \mathrm{m}$ resulting in $N_{\mathrm{p}}^1=456$, $N_{\mathrm{p}}^2=532$, and $N_{\mathrm{p}}^3=626$.  The conductivities of the inner and outer layers are $\sigma_1=1 \ \mathrm{S/m}$ and $\sigma_3=1 \ \mathrm{S/m}$, respectively. Thirteen sets of the conductivities of the middle layer are considered as $\sigma_2=\{50, 25, 10, 7.5, 5, 2.5, 1.25, 10, 0.75, 0.5, 0.25, 0.125, 0.1   \} \times 10^{-2} \ \mathrm{S/m}$, i.e., the  ratio of the conductivities between the brain/scalp layer and the skull layer is $\{2, 4, 10, 13.33, 20, 40, 80, 100, 133.33, 200, 400, 800, 1000 \}$, which covers the usual range of the conductivity ratio \cite{Goncalves}.

\begin{figure}[!t]
\centering
\includegraphics[width=0.9\columnwidth]{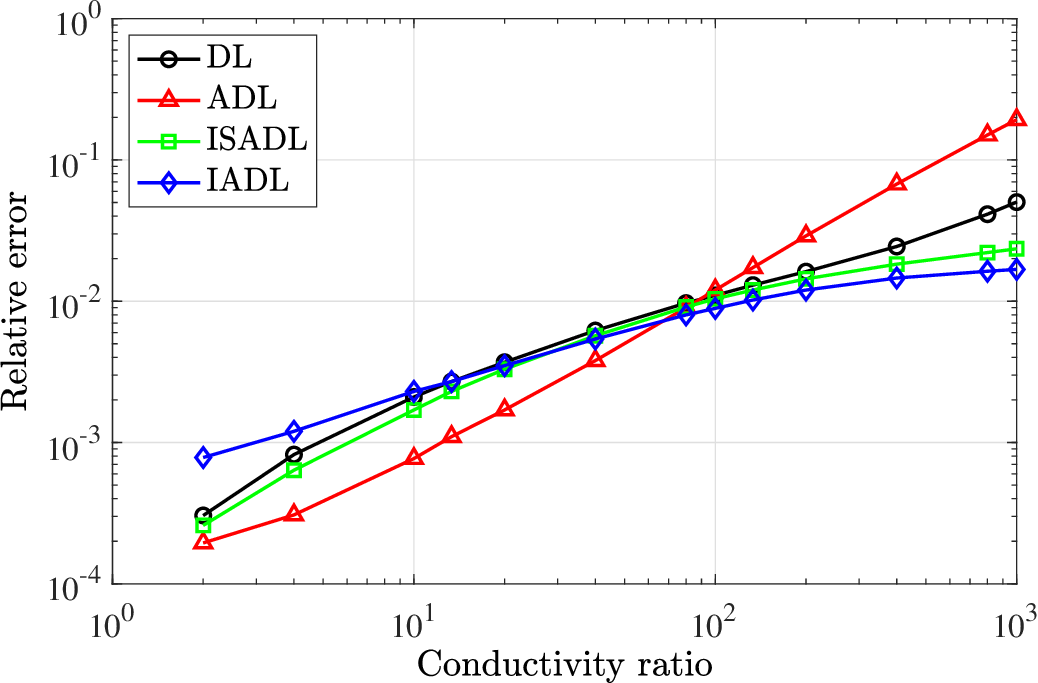}
\caption{ The relative error of the solution after using the proposed scheme  versus the conductivity ratio for the DL, ADL, ISADL, and IADL approaches.}
\label{fig_7}
\end{figure}

Fig. \ref{fig_7} compares the relative error of the solution after using the proposed scheme  versus the conductivity ratio for the DL, ADL, ISADL, and IADL approaches. This figure shows that, the relative error of the solution after using the proposed scheme for all  the approaches increases with the decrease of  $\sigma_2$, which is due to the isolated skull problem \cite{Hamalainen}. However, due to the improvement of the SIE formulation, the ISADL approach is more resistant to the increase of the conductivity ratio than the DL approach. Moreover, the IADL approach has the best accuracy among all the four approaches for large conductivity ratios, which matches with the behavior observed from Fig. 8 in \cite{Rahmouni}. Clearly, these results demonstrate the efficiency of the proposed scheme for discretizing SIEs for  all the four approaches.

\section{Conclusion}
\label{sec:conclusion}
In this work we propose a Nyström-based HO discretization scheme for SIEs for analyzing the EEG forward problem. The HO surface elements and the HO basis functions are used for the discretization of the interfaces of the head volume and the unknowns on the elements, respectively. Differently from existing works using the isoparametric HO discretization scheme, this work employs interpolation points different from the mesh nodes, which allows for the flexible manipulation of the order of the basis functions without the requirement of the regeneration
of the mesh of the interfaces. In addition, since the interpolation points are the same as the Gaussian quadrature points on the elements, the numerical computation of the surface integrals is simplified for the far-interaction case. Several numerical experiments are presented to  demonstrate the accuracy, flexibility, and efficiency
of the proposed discretization scheme for the DL, ADL, ISADL, and IADL approaches for analyzing the EEG forward problem.

The development of the Nyström-based HO discretization scheme for the symmetric SIE \cite{Kybic} for analyzing the EEG forward problem is undergoing, which needs the treatment of the hyper-singularity for the self-interaction case, and it will be presented in our future work.







\end{document}